\documentclass[11pt]{article}
\usepackage{graphicx}

\usepackage[letterpaper,width=135mm,height=203mm]{geometry}

\newenvironment{Proof}{\textsc{Proof.}}{\hspace{8mm}\QEDmark\smallskip}
\usepackage{latexsym}
\newcommand{\QEDmark}{\mbox{$\Box$}}

\newtheorem{theorem}{Theorem}
\newtheorem{lemma}[theorem]{Lemma}

\newcommand{\dd}[1]{\textbf{\textit{#1}}}

\usepackage{amsmath}

\newcommand{\curlyF}{\mathcal{F}}

\newcommand{\curlyT}{\mathcal{T}}

\newcommand{\aiso}{\zeta_k}
\newcommand{\iaiso}{\zeta_k^i}

\begin{document}


\begin{center}
{\Large \textbf{All-$k$-Isolation in Trees}} \\[3mm]
Geoffrey Boyer, Garrett C. Farrell, Wayne Goddard \\[1mm]
School of Mathematical and Statistical Sciences \\
Clemson University 
\end{center}

\begin{abstract}
We define an all-$k$-isolating set of a graph to be a set $S$ of vertices such that, if one removes $S$ and all its neighbors, then no component
in what remains has order $k$ or more. The case $k=1$ corresponds to a dominating set
and the case $k=2$ corresponds to what Caro and Hansberg called an isolating set.
We show that every tree of order $n \neq k$ contains an all-$k$-isolating set $S$ of size at most $n/(k+1)$, and moreover,
the set $S$ can be chosen to be an independent set. This extends previous bounds on variations of isolation,
while improving a result of Luttrell et al., who called the associated parameter the $k$-neighbor component order connectivity. 
We also characterize the trees where this bound is achieved. Further, we show that for~$k\le 5$,
apart from one exception every tree with $n\neq k$ contains $k+1$ disjoint independent all-$k$-isolating sets.
\end{abstract}

\section{Introduction}

In~\cite{CH}  Caro and Hansberg introduced a generalization of domination. Specifically, for a graph $G$, 
they measured the impact of a set $S \subseteq V(G)$ by what 
remains if $S$ and its neighbors (collectively denoted $N[S]$) is removed. For any forbidden set $\curlyF$ of graphs, a set $S$ is \dd{$\curlyF$-isolating}
if $G-N[S]$ contains no (not necessarily induced) copy of a graph in~$\curlyF$. If $\curlyF$ contains only one graph, say $F$, then 
we write $F$-isolating, and the $F$-isolation number is the minimum size of an $F$-isolating set. Domination corresponds to the case that $F = K_1$, 
and what is now just isolation corresponds to the case that $F = K_2$, also known as vertex-edge domination; see e.g.~\cite{Peters,LHHF}.
The original paper and other authors have considered
other possibilities for the set $\curlyF$, including stars, small paths \cite{CLWX,CXpath,CZZ,ZW}, cliques~\cite{BFK}, and all cycles \cite{BorgCycle,ZWcycle}.

In this paper we define an \emph{all-$k$-isolating  set} of $G$ to be a set $S$ of vertices such that $G-N[S]$ has no component of order at least~$k$. 
An all-$k$-isolating  set is thus equivalent to $\curlyF$-isolation when $\curlyF$ is all connected graphs of order $k$
or equivalently all trees of order $k$.  
For a graph $G$, we define its \emph{all-$k$-isolation  number} to be the minimum size of an all-$k$-isolating  set and denote it $\aiso(G)$. So domination
is all-$1$-isolating, ordinary isolation is all-$2$-isolating, and
all-$3$-isolating is equivalent to $P_3$-isolation or what is called $1$-isolation in~\cite{CH}.

We note that all-$k$-isolation is equivalent to ``$k$-neighbor component order connectivity'' introduced by Luttrell et al.~\cite{Luttrell12}, and 
an algorithm for trees was recently provided by Luttrell~\cite{Luttrell24}. But we propose the new terminology to link it to the isolation literature.
We further remark that the concept of removing the closed neighborhood of a set of vertices and caring about the maximum component
order that remains was also studied by Cozzens and Wu~\cite{CW} as the \emph{vertex neighborhood-integrity} of a graph.

At the same time, there are other ways to generalize domination.
For example, adding the requirement to a dominating set that it be independent produces the independent domination number. This parameter
has been extensively investigated. For example, while every graph without isolates has domination number at most $n/2$,
the independent domination number can be as large as $n - O(\sqrt{n})$, though for trees it is bounded by $n/2$. 
We~\cite{BGindep} showed a similar contrast between
isolation number and independent isolation number---while Caro and Hansberg proved that the ordinary version is at most $n/3$
except for $K_2$ and $C_5$, we observed that the independent isolation number can be asymptotically $n/2 - O(\sqrt{n})$, though
for trees (other than $K_2$) it is at most $n/3$. 

In this paper we extend these results for trees.
In this regard, Zhang and Wu~\cite{ZW} showed that for every connected graph of order $n\ge 4$ except $C_6$, it holds that
the $P_3$-isolation number is at most $2n/7$. Furthermore, they showed that there are infinitely many
graphs where the bound is attained, and the full extremal family is given in~\cite{CLWX,CZZ}.
For stars in general, Caro and Hansberg~\cite{CH} showed that the $K_{1,k-1}$-isolation number of a graph is at most $n/k$.
And if $G$ is a tree other than $K_{1,k-1}$ itself, they improved the upper bound to $n/(k+1)$.
(In their notation $\iota_k$ means $K_{1,k+1}$-isolation.) 
We continue
the investigation of these ideas in tree, in particular when there is the added requirement that the set $S$ be an independent set.
We denote by $\iaiso(G)$ the independent version, meaning where the all-$k$-isolating set is required to be an independent set.

Apart from bounds, there have also been results on colorings or partitions of graphs.
For example, it is well-known that a graph without isolates has two disjoint dominating sets. In~\cite{BGfirst} we showed that any connected
graph other than $K_2$ and $C_5$ has a partition into three disjoint isolating sets, while in~\cite{BGindep} we showed that  if the graph is bipartite
this partition
can be chosen to be a proper coloring. In particular, every tree has a partition into three disjoint independent
isolating sets.

We proceed here as follows. In Section~\ref{s:tree} we consider trees $T$ and provide a sharp upper bound of $n/(k+1)$ for $\aiso(T)$ and $\iaiso(T)$.
In Section~\ref{s:colorTree} we show that for $k\le 5$ a stronger result is possible, in that the tree can be partitioned into $k+1$ independent 
all-$k$-isolating sets (but this fails for $k\ge 7$). In Section~\ref{s:end} we provide some concluding remarks.

\section{An Upper Bound for Trees} \label{s:tree}

In this section we extend the upper bound $n/(k+1)$ of Caro and Hansberg~\cite{CH} on the
$K_{1,k-1}$-isolation number of a  tree of order $n$ in a couple of ways. First the bound for $K_{1,k-1}$-isolation
holds for any forbidden tree $F$ of order $k$; and 
second the bound holds for the independent version. The upper bound is also best possible for any $F$. 

It is to be noted that Luttrell et al.~\cite{Luttrell12,LuttrellThesis} provide a construction and an upper bound
for, in their terminology,  the $k$-neighbor component order connectivity of a tree. Moreover they claim the construction and upper bound match. 
However, Observation~4.5 of \cite{Luttrell12} is incorrect. In fact, the upper bound they obtain is $\aiso(T) \le \lfloor n / k \rfloor$
while the construction provides trees where $\aiso(T) = \lfloor n / (k+1) \rfloor$. We confirm here that the construction is 
optimal. We also determine the full family of extremal graphs in the case that $n$ is a multiple of $k+1$.

\begin{theorem} \label{t:iaisoTree}
    If $T$ is a tree of order $n\neq k$, then $\aiso(T) \leq \iaiso(T) \leq \frac{n}{k+1}$.
\end{theorem}
\begin{Proof}
It is known that the independent domination number of a tree is at most half its order. So assume $k\ge 2$.
For fixed $k$ we proceed by induction on $n$.    It is trivial that the result holds for $n<k$, so assume $n>k$.

For an independent set $J$ of vertices, let $S_k(J)$ denote the set of
vertices of $T-N[J]$ that are in components of order less than~$k$. 
The key observation is that:
\begin{quote}
If $t_J = |N[J] \cup S_k(J)| \ge (k+1)|J|$ and $T-N[J]$ has no component of order exactly $k$,
then one can induct on (the components of) $T-N[J]-S_k(J)$ and we are done.
\end{quote}
This follows since one can extend the independent
all-$k$-isolating set of $T-N[J]-S_k(J)$  by adding $J$ to form an independent all-$k$-isolating
set of $T$ with size at most $ t_J/(k+1) + (n-t_J)/(k+1) = n/(k+1)$.

Recall that a \dd{centroid} of a tree is a vertex $v$ such that the maximum component order of $T-v$ is as small
as possible. It is known that that maximum component order for centroid deletion is at most $n/2$. 
If $T$ has order at most $2k+1$, then let $v$ be the centroid. The set $\{v\}$ 
is thus an all-$k$-isolating set of the requisite size.
So assume $n \ge 2k+2$.

Consider any component that would be produced by removing one edge $e$. 
We say such a component is a \dd{$k$-branch} if it has order equal to $k$. 
Further, we say the end of $e$ in the $k$-branch is a \dd{branch} vertex
and the other end of $e$ is a \dd{core} vertex. (Note that two branch vertices cannot be adjacent by the assumption on the order.) 

\textbf{Case 1:} \textit{There is no $k$-branch.}
Then in particular, removing any connected subgraph of $T$ cannot leave a component of order exactly $k$. 
Define an edge to be \dd{large} if both components on its removal would have order at least $k+1$.

\textit{Case 1a: There is no large edge.} 
Then let $v$ be the centroid of $T$. 
For each neighbor $w$ of $v$, the subtree of $T-vw$ containing $w$ has order 
at most that of the subtree containing $v$. Since the edge joining $v$ and $w$ is not large,
it follows that every component of $T-v$ has order at most $k$, and thus
every component of $T-N[v]$ has order less than $k$; that is,
 $S_k (\{v\})$ is the whole of $T-N[v]$. Thus $\{ v \}$ is an all-$k$-isolating set.

\textit{Case 1b: There is a large edge.} 
Then consider a longest path in the subtree of~$T$ induced by the large edges. 
Let $v$ be an end
of that path, say with large edge going to neighbor $w$. We claim $|N[v] \cup S_k (\{v\})| > k+1$. 
For, since $v$ is incident with only one large edge, the set $ N[v] \cup S_k (\{v\})$ includes
all vertices in the component of $T-vw$ containing $v$, as well as $w$, which is enough. Further,
$T - N[v]$ has no component of order $k$. Hence one can induct.

\textbf{Case 2:} \textit{There exists a $k$-branch.}
 
\textit{Case 2a: There exists a core vertex $v$ that is adjacent to exactly one $k$-branch.}
Say $x$ is the branch neighbor of $v$. Then $|N[x] \cup S_k (\{x\})| \ge k+1$, and there is no component of order $k$
in $T-N[x]$, and so one can induct by the key observation.

\textit{Case 2b: There is no core vertex adjacent to exactly one $k$-branch.}
Consider a maximal subtree $C_T$ of $T$ all of whose vertices are core. 
Consider the bipartition $(B,R)$ of $C_T$ (where we assume $R=\phi$ if $C_T$ is $K_1$). 
See Figure~\ref{f:fBR}.
Let $b_k$ 
and $r_k$ denote the number of $k$-branches incident with
vertices of $B$ and $R$ respectively. Without loss of generality
assume $b_k \ge r_k$. Then form independent set $J$ by taking $B$ and 
the branch vertices (if any) adjacent to $R$. Note that $|J|=|B|+r_k$. 

    \begin{figure}[h]
        \centering
 \begin{center}
        \includegraphics[width=0.5\linewidth]{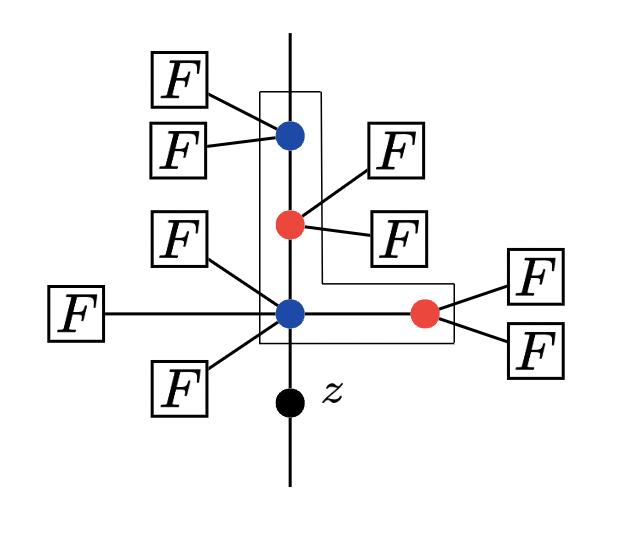}
 \end{center}
        \caption{The subgraph $C_T$ and its $k$-branches}
        \label{f:fBR}
    \end{figure}

Now, we claim $T-N[J]$ has no component of order $k$. For suppose there is such a component $F$.
Then $F$ must be a $k$-branch. Say $F$ is adjacent to vertex $z$ of $N[J]-J$. It follows that $z$ is core. 
Say $z$'s neighbor in $J$ is $w$. Then $w$ cannot be a branch vertex; so $w\in B$. But 
then this contradicts the maximality of the subtree $C_T$, since $z$ can be added to $C_T$. 
This proves the claim.

Further,  
every $k$-branch adjacent to $B \cup R$ contains a vertex of $N[J]$.
Thus $|N[J] \cup  S_k(J)| \ge  |B|+|R| + (b_k+r_k)k$. 
We claim $|N[J] \cup  S_k(J)|  \ge (k+1) |J|$.
For, consider the quantity
\[
     \psi = \left[ |B|+|R| + (b_k+r_k)k \right] - \left[ (k+1) ( |B| + r_k ) \right] =  
     -k |B| + |R|  + kb_k - r_k .
\]
As a function of $|R|$, the quantity $\psi$ is minimized at $|R|$ as small as possible, namely $|R|=0$. 
As a function of $r_k$, the quantity $\psi$ is minimized at $r_k$ as large as possible, namely $b_k$. That is,
\[ 
     \psi \ge -k|B| + (k-1)b_k .
\]
But $b_k \ge 2|B|$, since we are not in Case~2a. Hence $\psi \ge 0$ (using the fact that $k\ge 2$), as claimed.
Thus one can induct.
\end{Proof}

Let the  $F$-isolation number and independent $F$-isolation numbers of a tree~$T$ be denoted by
$\iota_F(T)$ and $\iota_F^i(T)$ respectively. Since the independent version is an upper bound for the standard version,
the following theorem follows immediately. 

\begin{theorem} \label{T:upperBound}
If $F$ is a tree of order $k$ and $T$ is a tree of order $n\neq k$, 
then $\iota_F (T ) \leq \iota^i_F ( T ) \leq \frac{n}{k+1}$.    
\end{theorem}

\begin{figure}[h]
\begin{center}
\includegraphics{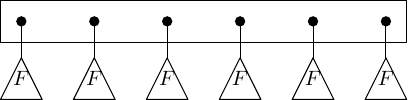}
\caption{A tree with $F$-isolation number $n/(k+1)$ for $F$ of order $k$}
\label{f:extreme}
\end{center}
\end{figure}

As suggested by the construction in~\cite{CH}, one can create an infinite family of extremal trees. 
Start with a base tree, say $T_0$. For each vertex $v$ of $T_0$, introduce
a copy of $F$, say $F_v$, and join $v$ to one vertex in $F_v$. See Figure~\ref{f:extreme}.
One can verify that this tree attains the bound and thus Theorem~\ref{t:iaisoTree} is sharp:
every $F$-isolating set contains at least one vertex from each $V(F_v) \cup \{v\}$.
Note that these extremal graphs also have the property that their independent $F$-isolation number is $\frac{n}{k+1}$.

More generally, one can  for each vertex
$v$ of the tree $T_0$, independently introduce
a copy of any tree of order $k$, say $F_v$, and join $v$ to one vertex in $F_v$. 
Let $\curlyT_k$ denote the resultant family of trees.

\begin{theorem}
For $k \ge 2$,
the following are equivalent for a tree $T$ of order $n$: \\
(i) $\aiso(T) = \frac{n}{k+1}$; \\
(ii) $\iaiso(T) = \frac{n}{k+1}$; \\
(iii) $T \in \curlyT_k$.
\end{theorem}
\begin{Proof}
It is immediate that for a tree $T$ in $\curlyT_k$, every all-$k$-isolating set contains at least one vertex from each $V(F_v) \cup \{v\}$.
Thus using Theorem~\ref{t:iaisoTree} it follows that condition~(iii) implies condition~(i), and condition~(i) implies condition~(ii).
To complete the proof of the theorem, it remains to show that if a tree has $\iaiso(T) = \frac{n}{k+1}$ then $T \in \curlyT_k$.

For fixed $k$ we proceed by induction on the order.
We need equality in the proof of Theorem~\ref{t:iaisoTree}. Equality in the base case is only attained for trees of order $k+1$.
Otherwise we need equality in the inductive step. 
It can be checked that the only place where the induction in that proof doesn't yield a smaller bound than 
the statement, is in Case 2a. This is immediate if $k>2$. 
If $k=2$, in Case~2b while we proved only that $\psi \ge 0$, it is immediate that 
achieving $\psi = 0$ is impossible since it requires $|R|=0$ and $r_k=b_k>0$, which is a contradiction.
Thus we need only worry about Case 2a. That is, $T$ must have a $k$-branch, say connected to core vertex $v$.
Let $U$ be the forest that remains when the $k$-branch and $v$ are removed.
For $T$ to have the maximum independent all-$k$-isolation number it must be that 
$\iaiso(U) = (n-k-1)/(k+1)$. Thus
by the induction hypothesis of this theorem, every component  of $U$ is in $\curlyT_k$. 

Note that vertex $v$ has exactly one neighbor in each component of $U$. We build an 
independent all-$k$-isolating set of $T$ as follows.
For each component $T'$ of $U$ that has order $s(k+1)$ for $s>1$, call  a vertex of the base tree of $T'$ 
a core vertex as before,
and call its unique non-core neighbor a branch vertex.
If $v$ is adjacent to none of the branch vertices in~$T'$, then let the ``representatives'' for $T'$ be its $s$ branch vertices.
If $v$ is adjacent to one of the branch vertices in~$T'$, then let the representatives for $T'$ be its branch vertices, 
except for the one adjacent to~$v$, along with the core vertex that shares a neighbor with $v$.
For each component $T'$ of $U$ that has order $k+1$, let one vertex nonadjacent to $v$ be its representative.
Let $D$ denote the set of all representatives. Note that $D \cup \{v \}$ dominates the core vertices
and is an independent all-$k$-isolating set of $T$
of size $n/(k+1)$.

Suppose vertex $v$ is adjacent to a non-core vertex $w$ in a component $T'$ of~$U$ of order more than $k+1$.
Then let $b$ be the core vertex of $T'$ chosen as representative. The set $( D - \{ b \} ) \cup \{ v \}$ 
is an independent all-$k$-isolating set of $T$, and thus $T$ is not extremal.
Similarly, suppose vertex $v$ is adjacent to a non-end vertex $w$ in a component $T'$ of order $k+1$.
Let $b$ denote the representative from that component. Again 
the set $( D - \{ b \} ) \cup \{ v \}$ is an independent all-$k$-isolating set.
That is, the core-vertex $v$ can only be adjacent to other core vertices or to an end-vertex
of each order-$(k+1)$-component.
Hence, $T$ is in $\curlyT_k$. \end{Proof}

Note that this characterization does not extend to independent domination (the case $k=1$)! For example, 
the double star has independent domination number half its order.

\section{Coloring Trees} \label{s:colorTree}

In this section 
we show that for small $k$ one can achieve something stronger than Theorem~\ref{t:iaisoTree}. 
It is immediate in a nontrivial tree that both partite sets are independent dominating sets.
In~\cite{BGindep} we showed that for trees with order $n\ge 3$ one can partition the 
vertex set into three independent isolating sets, or equivalently there is a proper $3$-coloring
where each color is an all-$2$-isolating set.
We generalize this result here for $k\le 5$.

To make the inductive step easier, we define an \dd{$\ell^*$-coloring} of a graph as a weak-partition of the vertex set into $\ell$ independent sets, that is, where
the sets are allowed to be empty but the coloring is proper.
Inspired by the definition in~\cite{JKOWdynamic}, we define an $\ell^*$-coloring of a graph $G$ as \dd{dynamic} if every vertex sees as
many of the colors as its degree allows: in other words, for all  vertices $v$, the number
of colors present in its (open) neighborhood is $\min( d_G(v), \ell-1 )$ where $d_G(v)$ is the degree of $v$.
The first observation is that if one has the desired coloring of a tree, then one may assume the
coloring is also dynamic. For any color $i$, we define an \dd{$i$-remnant} as any component
of the graph that remains when all vertices of color $i$ and their neighbors are removed. 

\begin{lemma} \label{l:dynamic}
For any fixed family $\curlyF$ of connected graphs,
if a tree $T$ has an $\ell^*$-coloring where every color is an $\curlyF$-isolating set of~$T$,
then there is such a coloring that is dynamic.
\end{lemma}
\begin{Proof}
Consider the $\ell^*$-coloring of $T$, and assume vertex $v$ sees less colors than its degree allows. 
Then there must be a color, say $1$, that is present more than once in $N[v]$ (necessarily on two of $v$'s neighbors), 
and a color, say $2$, that is missing from $N[v]$.
Let $w$ be a neighbor of $v$ of color~$1$.
Then within the subtree of $T-vw$ containing $w$, recolor every vertex with color $1$ with color $2$ and recolor every vertex with color 
$2$ (if any) with color $1$.

We claim that the new coloring still has every color class being $\curlyF$-isolating. Clearly, colors
other than $1$ and $2$ are unaffected. But since $v$ is now dominated by both
colors $1$ and $2$, any $1$- or $2$-remnant is completely contained within a component of $T-vw$ and thus
was there before. Thus there is no large remnant by the assumption on the initial coloring.
Furthermore, this interchange increases the number of colors that $v$ sees while
not changing the number of colors that $w$ sees. Thus, repeated application of this interchange yields a dynamic coloring
\end{Proof}

We will also repeatedly use the following strengthening of the dynamic idea. Specifically, for
some edge $xy$ we will color the two components $T_x$ and $T_y$ of $T-xy$ separately
and then merge the $\ell^*$-colorings. Before merging the colorings, we will re-number the colors 
so that after merger the coloring remains dynamic, and further, vertices $x$ and $y$ combined
see as many colors as possible. For definiteness, we will assume that $x$ 
has color $1$ and its $d_T(x)-1$ neighbors in $T_x$ have colors $2$ up to $\min( d_T (x) , \ell )$, while
vertex $y$ has color $\ell$ and its $d_T(y)-1$ neighbors in $T_y$ have colors $\ell-1$ down to 
$\max ( \ell - d_T( y )+1 , 1)$.
We say that we \dd{align} the colorings of $T_x$ and $T_y$. 

\begin{theorem} \label{t:fourColor}
Except for $P_3$, every tree $T$
has a $4^*$-coloring such that every color is all-$3$-isolating.
\end{theorem}
\begin{Proof}
The proof is by induction on the order $n$ of $T$.
The fact is immediate for $n\le 2$. So assume $n\ge 4$. 
Thus there exists an edge $e=xy$ such that $d_T(x)+d_T(y)\ge 4$.

Let $T_x$ and $T_y$ denote the components of $T-e$ containing $x$ and $y$ respectively. 
If $T_x$ is not $P_3$, then apply the induction hypothesis to it. By Lemma~\ref{l:dynamic} we may assume
the resulting coloring is dynamic. 
If $T_x$ is $P_3$, then give vertex $x$ color~$1$, (one of) its neighbor color~$2$, and the third vertex color $3$.
Color $T_y$ similarly, except reversing the ordering of the colors. 

Consider the combined aligned coloring in $T$. Note that this coloring is proper.
Vertices $x$ and $y$ combined see all $4$ colors, by the degree constraint. 
Thus the edge~$e$ is always dominated. Hence each $i$-remnant
must be completely contained within either $T_x$ or $T_y$. If $T_x$ is not~$P_3$, then by the inductive hypothesis, 
this remnant has order
at most two.   And if $T_x$ is~$P_3$, note that $y$ has the color missing on $T_x$;
thus at least one vertex of $T_x$ is dominated by each color and so what remains of $T_x$ has order at most two. 
An identical argument holds for $T_y$.
Hence the combined aligned coloring of $T$ has the desired property.
\end{Proof}

Note that Theorem~\ref{t:fourColor} is not true for bipartite graphs in general. For example, the $6$-cycle does not have such a coloring.


\begin{theorem} \label{t:fiveColor}
Except for trees of order $4$,
a tree $T$ has a $5^*$-coloring such that every color is an all-$4$-isolating set.
\end{theorem}
\begin{Proof}
We prove the statement by induction on the order of $T$. 
This is true for trees of order at most $3$.
So assume $T$ has order at least $5$. If $T$ is a path, then the natural sequential coloring suffices. So 
assume $T$ has a vertex $x$ of degree at least~$3$. If $T$ is a star, 
then any coloring using all five colors suffices. 
So assume $x$ has a neighbor $y$ of degree at least $2$.

Consider removing the edge $xy$, say yielding components $T_x$ containing $x$ and
$T_y$ containing $y$. We say a component is \dd{bad} if it is has order $4$; otherwise it is \dd{ordinary}.
By the induction hypothesis and Lemma~\ref{l:dynamic}, if $T_x$ is ordinary one can $5^*$-color it with a dynamic coloring
such that each color is an all $4$-isolating
set, and choose that vertex $x$ gets color $1$ and its neighbors in $T_x$ get colors $2$ up to $\min( d_T(x),5)$.
If $T_x$ is bad one can use colors $1$ through $4$ on it such that $x$ gets  color $1$ and its neighbors in $T_x$ get colors $2$ up to $d_T(x)$.
As usual we color $T_y$ similarly, except the ordering of the colors is reversed.

Now consider this coloring in $T$. 
Vertices $x$ and $y$ combined see all five colors; hence any $i$-remnant
must be completely contained within either $T_x$ or $T_y$.
If $T_x$ is ordinary, all its $i$-remnants are small.
If $T_x$ is bad, then its missing color, namely $5$, is on~$y$, 
and removing that color and its neighborhood leaves $T_x$ as a component of order at most $3$.
Similarly with $T_y$. Thus we have the desired coloring. 
\end{Proof}

Define $O_7$ as the tree of order $7$ obtained from $K_{1,3}$ by subdividing each edge once.

\begin{theorem} \label{t:sixColor}
Except for $O_7$ and trees of order $5$,
a tree $T$ has a $6^*$-coloring such that every color forms an all-$5$-isolating set.
\end{theorem}
\begin{Proof}
We prove the statement by induction on the order of $T$. 
This is true for trees of order at most $4$, by any $6^*$-coloring.
So assume the tree $T$ has order at least $6$. 
Consider the maximum degree-sum of adjacent vertices.\medskip

\textbf{Case 1:} \textit{There exist adjacent vertices $x$ and $y$ such that $d_T(x)+d_T(y) \ge 6$.}\\
Consider removing the edge $xy$, say yielding components $T_x$ containing $x$ and
$T_y$ containing $y$. We say a component is \dd{bad} if it is  $O_7$ or has order $5$; otherwise it is \dd{ordinary}.
We will color $T_x$ such that:
\begin{quote} \itshape
(a) vertex $x$ gets color $1$ and its neighbors in $T_x$  get colors $2$ up to $\min( d_T(x),6)$; and  \\[1mm]
(b) if $T_x$ is ordinary then every color is all-$5$-isolating; while 
if $T_x$ is bad then colors $1$ through $5$ are all-$5$-isolating, and the $6$-remnant
contains~$x$ and has order exactly $5$.
\end{quote}

If $T_x$ is ordinary, then by the induction hypothesis,  
one can $6^*$-color it such that each color is an all-$5$-isolating set of $T_x$.
By Lemma~\ref{l:dynamic},  we may assume the coloring is dynamic
and hence ensure the colors of $x$ and its neighbors are as specified.
If $T_x$ has order $5$, then color it with colors $1$ through~$5$, choosing the color of $x$ to be $1$ and its neighbor(s) 
to satisfy the condition.
If $T_x$ is $O_7$, then possible colorings of $O_7$ for each position of $x$ are shown
in Figure~\ref{f:o7}. Thereafter, color $T_y$ in a similar fashion, except that the ordering of the colors is reversed.

\begin{figure}[h]
\centerline{ \includegraphics[scale=1.1]{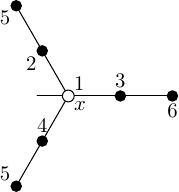}\qquad
\includegraphics[scale=1.1]{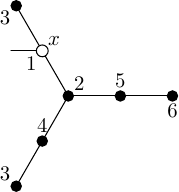}\qquad
\includegraphics[scale=1.1]{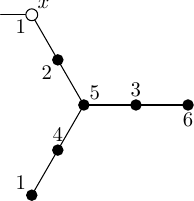}}
\caption{Possible colorings of $O_7$}
\label{f:o7}
\end{figure}

Now consider this coloring in $T$. 
By their degrees, vertices $x$ and $y$ combined see all six colors, and hence
if we remove from $T$ any color class and its neighborhood,
the edge $xy$ is destroyed. So it suffices to show that there is no remnant of order at least $5$ completely
contained within what remains of $T_x$ (since the same argument applies to $T_y$).
If $T_x$ is ordinary, then the inductive hypothesis says that this is true.
If $T_x$ is bad, then the coloring thereof ensures that the only possible bad case is a $6$-remnant. 
But such a subtree contains $x$, and color $6$ appears on~$y$, so $x$ is 
dominated by that color, and the $6$-remnant has order less than $5$.
Since each color is an all-$5$-isolating set,
we have the desired coloring.\medskip


\textbf{Case 2:} \textit{For all adjacent vertices $x$ and $y$ it holds that $d_T(x)+d_T(y) \le 5$.}\\
If $T$ is a path, then the natural sequential coloring suffices (since $T$ has order at least $6$). 
It cannot be the case that $d_T(x)=4$ because $T$ being $K_{1,4}$ is excluded.
So we may assume there is a vertex $x$ of degree
$3$ whose neighbors $y_1,y_2,y_3$ each have degree at most $2$. 

Define $T_i$ as the subtree induced by $x$ and $y_i$ and all vertices whose
path to $x$ contains $y_i$. Note that $x$ is an end-vertex in each~$T_i$.
Number the three subtrees such that: if there is a subtree isomorphic to $P_3$, then $T_1$ is such a tree.
In each subtree, if vertex $y_i$ has another neighbor, then call it $z_i$.

For each subtree $T_i$, we will give it a coloring where: 
\begin{quote} \itshape
(a) Vertex $x$ has color $1$ and vertex $y_i$ has color~$i+1$. 
If $z_1$ exists, it receives color $5$; while if $z_2$ and/or $z_3$ exist, they receive
color~$6$. And for all $i$, if $z_i$ has another neighbor,
then one such neighbor gets the color from $\{5,6\}$ that $z_i$ doesn't get. \\[1mm]
(b) if $T_i$ is ordinary then every color is all-$5$-isolating; while 
if $T_i$ is bad then five of the colors are all-$5$-isolating, and the bad color is not color~$1$.
\end{quote}
If $T_i$ is ordinary, then we apply the inductive hypothesis to it. By Lemma~\ref{l:dynamic} one may assume $x$ and $z_i$ (assuming it exists)
receive different colors. If $z_i$ has a neighbor $w$ with the same color as $x$, then one can recolor $x$ any color except its neighbor
and still have a valid coloring,
since any subtree of order $5$ in $T_i$ containing~$x$ necessarily contains $z_i$.
Thus we may assume vertices $x$, $y_i$, $z_i$, and~$w$ receive different colors, and so
we may renumber the colors to satisfy the above constraints.
If $T_i$ is a component of order $5$, then
we use five colors as per the above constraints.
If $T_i$ is an $O_7$, then we use all six colors as per the above constraints,
with the repeat color being color $1$ (on $x$ and another end-vertex). 
Figure~\ref{f:part} shows an example of part of a possible resultant coloring of $T$.

\begin{figure}
\centerline{\includegraphics[scale=1.1]{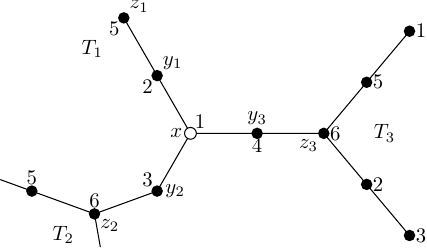}}
\caption{Part of coloring of $T$ around vertex $x$}
\label{f:part}
\end{figure}

Now, consider this coloring in $T$. 
Consider first a color $C \in \{ 1, \ldots, 4 \}$. 
Vertex $x$ is dominated by $C$. So any $C$-remnant is totally contained within some $T_i-x$.
If a $T_i$ is ordinary, then immediately $C$ is an all-$5$-isolating set of $T_i$.
If a $T_i$ is a component of order $5$, then at least one vertex of~$T_i$ is dominated by~$C$ (namely~$x$)
and hence the $C$-remnant has order less than~$5$. If a $T_i$ is an $O_7$, 
then $C$ dominates $x$ and at least one vertex nonadjacent to $x$;
again the $C$-remnant has order less than $5$.
That is, in all cases the vertices of color $C$ form an all-$5$-isolating set of~$T$.

Consider second a color $C \in \{5,6\}$. 
If a  subtree $T_i$ is ordinary, then there cannot be a $C$-remnant of order $5$ or more
totally contained within $T_i$. 
If $T_i$ is bad, then the central vertex $z_i$ is dominated by $C$.
That is, there cannot be a $C$-remnant of order $5$ or more totally contained within
 $T_i$.

Finally, we need to check that for $C \in \{5,6\}$ there  is no large $C$-remnant that straddles multiple $T_i$'s. 
Let $U$ be a subtree of order $5$ containing $x$ and contained within some $C$-remnant.
Such a subtree must contain a vertex
at distance two from $x$, say $z_j$. Furthermore $z_j$ must have 
degree one in~$T$, else it is guaranteed to be dominated by color $C$.
In particular, $T_j$ is isomorphic to $P_3$.
Since in this case we are assuming the maximum degree
of $T$ is $3$, there are only two possibilities for the isomorphism class
of $U$: either $P_5$ or the tree $Q$ of order $5$ with maximum degree $3$.

Suppose  $U$ is isomorphic to $P_5$. Then it 
must contain a second vertex at distance two from 
$x$, say $z_k$, which again has degree one in $T$. Since $T$ is not~$O_7$,
there is at least one $T_i$ that is not a $P_3$. It follows that one
of $z_j$ and $z_k$ is~$z_1$. Thus one of these vertices has color~$5$ and one
has color $6$, a contradiction of the claim that neither is dominated by~$C$.

Suppose $U$ is isomorphic to $Q$.
Then note that  $Q$ must consist of vertex $x$, its three neighbors, and $z_j$.
Since $T$ has order more than five, another neighbor of~$x$, say $y_k$, has 
degree two in $T$. For this copy of $Q$ to survive, we thus need $z_k$ to
not be color $C$ and so it must
have the same color as $z_j$. Hence $\{j,k\} = \{2,3\}$. But since $T_j$ is $P_3$,
by the choice of numbering of the $T_i$, that means that $T_1$ 
is also isomorphic to $P_3$. Thus $z_1$ exists, has color $C$, and dominates $y_1$,
a contradiction. That is, it follows that color $C \in \{5,6\}$ is also all-$5$-isolating, as required.
 \end{Proof}

We note that the pattern of Theorems~\ref{t:fourColor}, \ref{t:fiveColor}, and \ref{t:sixColor}
does not hold in general. That is, one cannot always
find $k+1$ disjoint all-$k$-isolating sets in a tree. We make no attempt to optimize the bound
in the following lemma.

\begin{lemma}
For each $k\ge 7$, there exist arbitrary large trees where the maximum number of disjoint $P_k$-isolating
sets is at most $(3k+10)/4 < k+1$.
\end{lemma}
\begin{Proof}
Assume $k$ is odd. Then  
let $H_k$ be the tree of order $(3k-1)/2$ obtained from $K_{1,3}$ with center $v$ by subdividing each edge $(k-3)/2$ times.
Let $w$ be an end-vertex of $H_k$ and consider
any tree $T$ that is formed from a tree $T_0$ by adding a disjoint copy of $H_k$ and joining $w$ to one vertex 
of $T_0$. Consider a $P_k$-isolating set $S$ of $T$. It is immediate that $S$ contains at least one vertex of $H_k$.
The key observation is that unless $S$ contain one of the four vertices of~$N[v]$ 
or the neighbor of $w$ in $T_0$, it must contain at least
two vertices from $H_k$. Thus the maximum number of disjoint $P_k$-isolating sets of $T$ is at most $5 + (|H_k|-5)/2 = 
(3k+9)/4$.   

The case that $k$ is even is similar. We let $H_k$ be the tree of order $3k/2$ obtained from $K_{1,3}$ by subdividing two
of the edges $(k-2)/2$ times and subdividing one edge $(k-4)/2$ times. Again consider any tree $T$ formed
from a tree $T_0$ by adding a copy of $H_k$ and joining an end-vertex of $H_k$ to $T_0$. By similar reasoning, the maximum number of disjoint
$P_k$-isolating sets of $T$ is at most $5 + (|H_k|-5)/2 = 
(3k+10)/4$.
\end{Proof}

It is unclear what happens for $k=6$. On the other hand,
while the all-$k$-isolating coloring result does not extend to all $k$, we observe that for the star it does, thereby strengthening the original bound
of Caro and Hansberg~\cite{CH}:

\begin{theorem} \label{t:colorStar}
For $k\ge 2$, 
except for $K_{1,k-1}$ itself, every tree $T$
has a $(k+1)^*$-coloring such that every color is $K_{1,k-1}$-isolating.
\end{theorem}
\begin{Proof}
We re-use the proof of Theorem~\ref{t:fourColor}.
The proof is by induction on the order $n$ of $T$.
The fact is immediate for $n<k$. So assume $n > k$. 
If the graph has maximum degree less than $k-1$, then any coloring suffices. So assume there
is a vertex $x$ with degree at least $k-1$. If $T$ is a star, then any coloring that uses all $k+1$ colors suffices.
So assume $T$ is not a star. Then $x$ has a neighbor $y$ with $d_T(y)\ge 2$.

As usual, let $T_x$ and $T_y$ denote the components of $T-e$ containing $x$ and $y$ respectively. 
If $T_x$ is not $K_{1,k-1}$, then apply the induction hypothesis to it. 
If $T_x$ is $K_{1,k-1}$, then color $x$ with color $1$ and use colors $2$ up to $k$ for the other vertices. 
Color $T_y$ similarly, except reversing the ordering of the colors.

Consider the combined aligned coloring in $T$. 
Vertices $x$ and $y$ combined see all colors, by their degrees. 
Thus the edge~$e$ is always dominated and each $i$-remnant
must be completely contained within either $T_x$ or $T_y$. If $T_x$ is not~$K_{1,k-1}$, then by the inductive hypothesis, 
this remnant does not contain $K_{1,k-1}$.   And if $T_x$ is~$K_{1,k-1}$, note that $y$ has the color missing on $T_x$, namely $k+1$.
An identical argument holds for $T_y$.
Hence the combined aligned coloring of $T$ has the desired property.
\end{Proof}

\section{Concluding Remarks} \label{s:end}

While the upper bounds for the ordinary and independent versions are the same, we observe that for a particular
graph the parameters can be very different.

\begin{theorem}
For any family $\curlyF$ of connected graphs, the 
independent $\curlyF$-isolation and ordinary $\curlyF$-isolation numbers of a graph can be arbitrarily far apart.
\end{theorem} 
\begin{Proof}
Let $F$ be a graph in $\curlyF$ of smallest order, and designate some vertex $f$ of~$F$ as its root.
Then take $2b$ disjoint copies of $F$, and introduce a copy of $K_2$ with vertex
set $\{x, y\}$ where $x$ adjacent to the roots of $b$ of the copies and 
 $y$ is adjacent to the roots of the other $b$ copies. For $\curlyF$-isolation
one can take $\{x,y\}$; so the number is~$2$. 
For independent $\curlyF$-isolation, we take $x$ and then need one vertex from each copy of $F$
adjacent to $y$, and so the number is $b+1$.
\end{Proof}

It remains to determine exactly which trees have the property that there is always a suitable coloring.
And an obvious further question is to establish similar bounds and coloring results for general graphs. For example,
perhaps the bound of Theorem~\ref{t:iaisoTree} holds for all graphs with sufficiently large girth.

\end{document}